\documentclass[a4paper,12pt]{article}
\usepackage{amsmath}
\usepackage{amssymb}
\usepackage{latexsym}
\usepackage{lscape}
\usepackage{pdflscape}
\usepackage[dvipsnames]{xcolor}
\usepackage{amsfonts}
\usepackage{amsmath}
\usepackage{latexsym}
\usepackage{enumitem}
\usepackage{graphics}
\usepackage{bbm}
\usepackage{tensor}
\usepackage{pdflscape} 
\usepackage{hyperref}
\usepackage{fixltx2e}  
\newcommand{\be}{\begin{eqnarray}}
\newcommand{\ee}{\end{eqnarray}}
\newcommand{\ben}{\begin{equation*}}
\newcommand{\een}{\end{equation*}}
\newcommand{\bean}{\begin{eqnarray*}}
\newcommand{\eean}{\end{eqnarray*}}
\def\bga#1\ega{\begin{gather}#1\end{gather}}
\def\bal#1\eal{\begin{align}#1\end{align}}
\newcommand{\ph}{\phantom}
\newcommand{\bsub}{\begin{subequations}}
\newcommand{\esub}{\end{subequations}}

\newcommand{\disfrac}[1][2]{\displaystyle\frac}

\newcommand{\non}{\nonumber}
\newcommand{\bs}{\boldsymbol}

\newcommand{\sech}{\textrm{sech\,}}
\newcommand{\ima}{\mathbbmtt{i}}
\newcommand{\IR}{\mathbbmtt{R}}

\newcommand{\con}{\lrcorner}
\newcommand*\xbar[1]{%
  \hbox{%
    \vbox{%
      \hrule height 0.5pt 
      \kern0.5ex
      \hbox{%
        \kern-0.1em
        \ensuremath{#1}%
        \kern-0.1em
      }%
    }%
  }%
}
\numberwithin{equation}{section}

\begin{document}

\title{\textbf{Faithful representations of Lie algebras and Homogeneous Spaces}}
\vspace{1cm}
\author{\textbf{Petros A. Terzis}\thanks{pterzis@phys.uoa.gr}\\
{\it Nuclear and Particle Physics Section, Physics Department,}\\{\it University of Athens, GR 157--71 Athens, Hellas}}
\date{}
\maketitle
\begin{center}
\textit{}
\end{center}
\vspace{-1cm}
\abstract{\textit{The faithful representations of real Lie algebras $\mathfrak{g}$ with dimensions $n\leq 4$ acting on the dual space of the universal covering algebra of an ideal $\mathfrak{s}$ of $\mathfrak{g}$, are calculated. With the help of left and right translations of the Lie group, we define the homogeneous manifolds (Lie groups with trivial isotropy subgroup) of dimension $n$ and we calculate the Killing fields (generators of left translations), the invariant fields (generators of right translations) and the left invariant 1-forms. The presented method also allow us to calculate the composition function of the Lie group. All of the calculated quantities are derived by applying algebraic manipulations and function differentiation.}}

\section{Introduction}

The notion of Lie groups and Lie algebras it encouters in a very broad area of modern mathematical physics. Abstract Lie algebras are widely used in quantum mechanics, e.g. to determine the energy spectrum of a system; the harmonic oscillator being the most notable example. In general relativity abstract Lie algebras play an important role when someone is constructing/classifying spacetimes which admit a set of Killing vector fields; the Bianchi classification \cite{Bianchi} is a prototype for this approach, see \cite{Jantzen1:Bianchi}, \cite{Jantzen2:Bianchi} for an English translation. The classification of real Lie algebras is a quite active area of research. The finite dimensional real analytic Lie algebras of vector fields on $\mathbb R^2$ were completely classified in \cite{Kamran}. The classification of finite-dimensional subalgebras in polynomial Lie algebras of rank one was made in \cite{Arzhantsev}; while all the finite-dimensional Lie algebras that can be realized as subalgebras of the Lie algebra of planar vector fields with coefficients from the field of rational functions were found in \cite{Makedonskyi}. A subjective overview of transitive Lie algebras and the classification of the primitive ones was the subject of \cite{Draisma}.

Modern applications of Lie algebras can be found in string theories; a beautiful branch is the Poisson--Lie T--duality introduced by Klimcik and Severa \cite{KlimcikSevera1}, \cite{Klimcik2} see also \cite{Sfetsos:1997pi} and \cite{Mojaveri:4+1hsc}. In \cite{ChrisDim} the automorphism group of 4d homogeneous manifolds was used for their invariant characterization. On the other hand, realizations of Lie algebras with vector fields are used in the integration of systems of ordinary and/or partial differential equations \cite{OlverODE}, a viewpoint that was used by the inventor of the theory, S. Lie. The interactive use of both the abstract and the realization approach was exploited in \cite{ChrisTerzBIII}, \cite{ChrisTerzBIII:2}, \cite{ChrisTerzBVII}, \cite{ChrisTerzBI-V}, resulting in the complete integration of the vacuum Einstein's field equations for several Bianchi cosmologies. Recently Nucci and Leach used Lie groups in order to analyze the inconsistencies in various schemes of quantization \cite{LeachNucci:Lie} and Nucci proposed a new way of quantization \cite{Nucci:Noether} based on Noether symmetries.

In this paper we use Lie' s work in a twofold way. Firstly, we use Lie groups, as abstract ones in order to define the homogeneous manifolds by the action of the group on itself with the help of right and left translations. Secondly, we use the Lie algebras accosiated to the Lie groups in order to construct faithful representations of these algebras, and finally by exponentiate these representations we evaluate the finite elements of the corresponding Lie group. At this stage we are able to calculate the Killing fields, the invariant bases and the basis 1-forms for the homogeneous manifold in question. The great advantage of our approach is that is completely algebraical; we only need the structure constants of the Lie algebra in order to build up its faithful representation acting on the dual space of the universal covering algebra. The method is algorithmic so it can be easily realized as a procedure in a computer algebra package.

The paper is organized as follows: Section 2 is divided in two subsections; the first is devoted to fixing the notation and constructing the homogeneous manifolds; while in the second we present our method. In Section 3 we put the method into work for 2d, 3d and 4d real Lie algebras, thus finding the representations of these algebras along with the Killing fields, the invariant bases and the basis 1-forms for the homogeneous manifolds. Finally, in the appendix we present the way to calculate the matrix exponential $\exp A$ of a matrix $A$ based on the definition of functions of matrices.

\section{Homogeneous manifolds}
In this section we present some basic elements from the theory of Lie groups and Lie algebras, in order to set up the notation along with the need to make the presentation self contained; for further details we refer the reader to classic works of Jacobson \cite{JacobsonLie} and Humphreys \cite{HumphreyLie}.

\subsection{General considerations}

Let $U_e$ be an open neighborhood of the identity element $e$ of an $N$--dimensional real Lie group $G$ and $\psi$ a homeomorphism of $U_e$ into an open set $V\subset \mathbb{R}^N$. Each element of the group can be parametrized with an $N$--tuple $\alpha=(\alpha^1,\alpha^2,\dots,\alpha^N) \in V$, i.e. $\alpha=\psi(g)$, thus we indicate this fact by writing
\be
\forall g\in G,\,g=g_\alpha.
\ee
The group operation $\diamond\colon G\times G\mapsto G$ reads
\bal
\forall g_\alpha,\,g_\beta\in G,\, g_\alpha\diamond g_\beta=g_{\phi(\alpha,\beta)},
\eal
where the $N$-- dimensional composition function $\phi(\alpha,\beta)$ defines the operation $\diamond$ on the coordinates $\alpha,\beta$
\bal
\phi \colon V\times V\mapsto V,\, \phi(\alpha,\beta)=\alpha\diamond \beta= (\phi^1,\phi^2,\dots,\phi^N),
\eal
and obey the relation
\bal\label{φφ σχέση}
\phi\left(\phi(\alpha,\beta),\gamma\right)=\phi\left(\alpha,\phi(\beta,\gamma)\right).
\eal

The Lie group $G$ can be made, in an intrinsic way, to act as a \emph{Lie group of transformations} on itself with the help of left $L_g$ and right $R_g$ translations
\bal\label{left right}
L_g\colon G\mapsto G,\, L_{g_x} g_\alpha=g_x\diamond g_\alpha,\quad R_g\colon G\mapsto G,\, R_{g_x} g_\alpha=g_\alpha\diamond g_x.
\eal

When we apply the above maps to a group element $g_\alpha=g_\varepsilon$  near the neutral element $e=e_0$ of $G$, we have the Taylor expansions for each component $\phi^\sigma$ of the composition function $\phi$
\bal
\phi^\sigma(x,\varepsilon)=\phi^\sigma(x,0)+ \left.\frac{\partial\phi^\sigma(x,\beta)}{\partial\beta^\lambda}\right|_{\beta^\lambda=0}\varepsilon^\lambda\Rightarrow\phi^\sigma(x,\varepsilon)=x^\sigma+\varepsilon^\lambda\bs{\xi}_\lambda x^\sigma\\
\phi^\sigma(\varepsilon,x)=\phi^\sigma(0,x)+ \left.\frac{\partial\phi^\sigma(\beta,x)}{\partial\beta^\lambda}\right|_{\beta^\lambda=0}\varepsilon^\lambda\Rightarrow \phi^\sigma(\varepsilon,x)=x^\sigma+\varepsilon^\lambda\bs{\eta}_\lambda x^\sigma,
\eal
where the above equalities are understood modulo $O(|\varepsilon|)$, as $\varepsilon\to 0$ and the Einstein summation convention is thereafter assumed. Furthermore, the vector fields $\{\bs{\xi}_\lambda,\bs{\eta}_\lambda\}\in T_xG$, the tangent space of $G$ at $x$, are defined by
\bsub\label{gen lr}
\bal
\bs{\xi}_\lambda=\left.\frac{\partial\phi^\tau(x,\beta)}{\partial\beta^\lambda}\right|_{\beta^\lambda=0}\frac{\partial}{\partial x^\tau}=\tensor{\xi}{_\lambda^\tau}(x)\partial_\tau \quad \text{with} \quad \tensor{\xi}{_\lambda^\tau}(0)=\tensor{\delta}{_\lambda^\tau}\\
\bs{\eta}_\lambda=\left.\frac{\partial\phi^\tau(\beta,x)}{\partial\beta^\lambda}\right|_{\beta^\lambda=0}\frac{\partial}{\partial x^\tau}=\tensor{\eta}{_\lambda^\tau}(x)\partial_\tau \quad \text{with} \quad \tensor{\eta}{_\lambda^\tau}(0)=\tensor{\delta}{_\lambda^\tau},
\eal
\esub
and are the generators of the left and right translations respectively. The commutator between $\bs{\xi}_\kappa$ and $\bs{\eta}_\lambda$ can be easily computed with the aid of \eqref{φφ σχέση} for the infinitesimals elements $\alpha^\lambda=\beta^\lambda=\varepsilon^\lambda$
\bal\label{comξη}
\phi^\sigma\left(\phi^\mu(\varepsilon^\kappa,x^\kappa),\varepsilon^\kappa\right)&=\phi^\sigma\left(\varepsilon^\kappa,\phi^\mu(x^\kappa,\varepsilon^\kappa)\right)\Rightarrow\non\\
\phi^\sigma(x^\mu+\varepsilon^\lambda\bs{\eta}_\lambda x^\mu,\varepsilon^\kappa)&=\phi^\sigma(\varepsilon^\kappa,x^\mu+\varepsilon^\mu\bs{\xi}_\mu x^\kappa)\Rightarrow\non\\
x^\sigma+\varepsilon^\lambda\bs{\eta}_\lambda x^\sigma+\varepsilon^\tau\bs{\xi}_\tau\left(x^\sigma+\varepsilon^\lambda \bs{\eta}_\lambda x^\sigma\right)&=x^\sigma+\varepsilon^\mu\bs{\xi}_\mu x^\sigma+\varepsilon^\tau\bs{\eta}_\tau\left(x^\sigma+\varepsilon^\mu\bs{\xi}_\mu x^\sigma\right)\Rightarrow\non\\
\varepsilon^\lambda\varepsilon^\kappa\left(\bs{\xi}_\kappa\bs{\eta}_\lambda- \bs{\eta}_\lambda\bs{\xi}_\kappa\right)x^\sigma&=0 \quad \forall x^\sigma\Rightarrow\non\\
\left[\bs{\xi}_\kappa,\bs{\eta}_\lambda\right]&=0,
\eal
i.e. each generator of left translations commutes with each generator of right translations; a well-known fact.

From the fundamental theorems of Sophus Lie we know that the two sets of the generators $\left\{\bs{\xi}_\kappa\right\}$ and $\left\{\bs{\eta}_\kappa\right\}$ span two Lie algebras, which are isomorphic, i.e.
\bal\label{comξcomη}
\left[\bs{\xi}_\kappa,\bs{\xi}_\lambda\right]=\tensor{C}{^\mu_{\kappa\lambda}}\bs{\xi}_\mu, \quad
\left[\bs{\eta}_\kappa,\bs{\eta}_\lambda\right]=\tensor{D}{^\mu_{\kappa\lambda}}\bs{\eta}_\mu,
\eal
where $\tensor{C}{^\mu_{\kappa\lambda}},\tensor{D}{^\mu_{\kappa\lambda}}$ are the structure constants for the corresponding algebra. We can calculate the values of $\tensor{D}{^\mu_{\kappa\lambda}}$ if we know the values of $\tensor{C}{^\mu_{\kappa\lambda}}$ as follows.

The vector fields $\bs{\xi}_\kappa$ span a basis of $T_xG$ thus the $\bs{\eta}_\lambda$ must be a linear combination of them with $x$-dependent coefficients
\bal\label{ηfromξ}
\bs{\eta}_\lambda=\tensor{c}{_\lambda^\tau}(x)\bs{\xi}_\tau.
\eal
Evaluating the above equation at the identity element $e=e_0$ and using \eqref{gen lr} we have
\bal\label{c0}
\tensor{\eta}{_\lambda^\mu}(0)=\tensor{c}{_\lambda^\tau}(0)\tensor{\bs{\xi}}{_\tau^\mu}(0)\Rightarrow \tensor{\delta}{_\lambda^\mu}=\tensor{c}{_\lambda^\tau}(0)\tensor{\delta}{_\tau^\mu}\Rightarrow \tensor{c}{_\lambda^\tau}(0)=\tensor{\delta}{_\lambda^\tau}.
\eal
Inserting now \eqref{ηfromξ} into the first of \eqref{comξη} and evaluating at $e=e_0$ we arrive at
\bal
\left[\bs{\xi}_\kappa,\tensor{c}{_\lambda^\tau}\right]\bs{\xi}_\tau+\tensor{c}{_\lambda^\tau}\left[\bs{\xi}_\kappa,\bs{\xi}_\tau\right]&=0\Rightarrow\non\\
\left[\bs{\xi}_\kappa,\tensor{c}{_\lambda^\tau}\right]\bs{\xi}_\tau+\tensor{c}{_\lambda^\tau}\tensor{C}{^\mu_{\kappa\tau}}\bs{\xi}_\mu&=0\Rightarrow\non\\
\left.
\left[\bs{\xi}_\kappa,\tensor{c}{_\lambda^\tau}\right]\right|_{x^\mu=0}\bs{\xi}_\tau(0) +\tensor{c}{_\lambda^\tau}(0)\tensor{C}{^\mu_{\kappa\tau}}\bs{\xi}_\mu(0)&=0\Rightarrow\non\\
\left.\left[\bs{\xi}_\kappa,\tensor{c}{_\lambda^\tau}\right]\right|_{x^\mu=0}&= -\tensor{C}{^\tau_{\kappa\lambda}},
\eal
where in the last line we used \eqref{c0} and \eqref{gen lr}. Finally, inserting  \eqref{ηfromξ} into the second of \eqref{comξη} and evaluating again at $e=e_0$ we end up with
\bal
\left[\bs{\eta}_\kappa,\bs{\eta}_\lambda\right]&=\tensor{c}{_\kappa^\mu}\left[\bs{\xi}_\mu,\tensor{c}{_\lambda^\rho}\right]\bs{\xi}_\rho+\tensor{c}{_\kappa^\mu}\tensor{c}{_\lambda^\rho}\left[\bs{\xi}_\mu,\bs{\xi}_\rho\right]+\tensor{c}{_\lambda^\rho}\left[\tensor{c}{_\kappa^\mu},\bs{\xi}_\rho\right]\bs{\xi}_\mu\Rightarrow\non\\
\left.\left[\bs{\eta}_\kappa,\bs{\eta}_\lambda\right]\right|_{x^\mu=0}&=\tensor{\delta}{_\kappa^\mu}\left(-\tensor{C}{^\rho_{\mu\lambda}}\right)\bs{\xi}_\rho(0)+\tensor{\delta}{_\kappa^\mu}\tensor{\delta}{_\lambda^\rho}\tensor{C}{^\tau_{\mu\rho}}\bs{\xi}_\tau(0)+\tensor{\delta}{_\lambda^\rho}\left(\tensor{C}{^\mu_{\rho\kappa}}\right)\bs{\xi}_\mu(0)\Rightarrow\non\\
\tensor{D}{^\mu_{\kappa\lambda}}\bs{\eta}_\mu(0)&=-\tensor{C}{^\rho_{\kappa\lambda}}\bs{\xi}_\rho(0)+\tensor{C}{^\rho_{\kappa\lambda}}\bs{\xi}_\rho(0)-\tensor{C}{^\mu_{\kappa\lambda}}\bs{\xi}_\mu(0),
\eal
resulting into
\bal\label{struD}
\tensor{D}{^\mu_{\kappa\lambda}}&=-\tensor{C}{^\mu_{\kappa\lambda}}.
\eal

Let $\{\bs{\sigma}^\lambda\}$ be the dual basis of $\{\bs{\eta}_\kappa\}$, i.e.
\bal\label{contr}
\langle\bs{\eta}_\kappa,\bs{\sigma}^\lambda\rangle=\bs{\eta}_\kappa\con\bs{\sigma}^\lambda= i_{\bs{\eta}_\kappa}\bs{\sigma}^\lambda=\tensor{\delta}{_\kappa^\lambda}.
\eal

For the contraction operator $\con$ with the vector fields $\bs{\eta}_\kappa$ and the $(p+1)$-forms $\bs{w},\,\bs{z}$ we have the following
\bal\label{w=z}
\bs{\eta}_\kappa\con\bs{w}=\bs{\eta}_\kappa\con\bs{z}\Rightarrow\bs{w}=\bs{z}.
\eal
To see that, write in components, $\tensor{\eta}{_\kappa^\lambda}w_{\lambda j_1\dots j_p}=\tensor{\eta}{_\kappa^\lambda}z_{\lambda j_1\dots j_p}$ and multiply each side with $\tensor{\sigma}{_\alpha^\kappa}$ in order to arrive to \eqref{w=z}.

The Lie drag of $\bs{\sigma}^\lambda$ with respect to $\bs{\xi}_\kappa$ and $\bs{\eta}_\kappa$ can easily be computed with the help of \eqref{contr}, \eqref{comξcomη} and \eqref{comξη}
\begin{gather}
\pounds_{\bs{\xi}_\kappa}\bs{\eta}_\lambda\con\bs{\sigma}^\rho=0\Rightarrow \left(\pounds_{\bs{\xi}_\kappa}\bs{\eta}_\lambda\right)\con\bs{\sigma}^\rho+ \bs{\eta}_\lambda\con\left(\pounds_{\bs{\xi}_\kappa}\bs{\sigma}^\rho\right)=0 \overset{\eqref{comξη}}{\Rightarrow} \bs{\eta}_\lambda\con\left(\pounds_{\bs{\xi}_\kappa}\bs{\sigma}^\rho\right)=0
\overset{\eqref{w=z}}{\Rightarrow}\non\\
\pounds_{\bs{\xi}_\kappa}\bs{\sigma}^\rho=0\label{lieξσ},
\end{gather}
and
\begin{gather}
\pounds_{\bs{\eta}_\kappa}\bs{\eta}_\lambda\con\bs{\sigma}^\rho=0\Rightarrow \left(\pounds_{\bs{\eta}_\kappa}\bs{\eta}_\lambda\right)\con\bs{\sigma}^\rho+ \bs{\eta}_\lambda\con\left(\pounds_{\bs{\eta}_\kappa}\bs{\sigma}^\rho\right)=0 \overset{\eqref{comξcomη}}{\Rightarrow} \non\\
-\tensor{C}{^\mu_{\kappa\lambda}}\bs{\eta}_\mu\con\bs{\sigma}^\rho+ \bs{\eta}_\lambda\con\left(\pounds_{\bs{\eta}_\kappa}\bs{\sigma}^\rho\right)=0
\overset{\eqref{contr}}{\Rightarrow}-\tensor{C}{^\rho_{\kappa\lambda}}+ \bs{\eta}_\lambda\con\left(\pounds_{\bs{\eta}_\kappa}\bs{\sigma}^\rho\right)=0 \overset{\eqref{contr}}{\Rightarrow}\non\\
-\tensor{C}{^\rho_{\kappa\tau}}\bs{\eta}_\lambda\con\bs{\sigma}^\tau+ \bs{\eta}_\lambda\con\left(\pounds_{\bs{\eta}_\kappa}\bs{\sigma}^\rho\right)=0 \overset{\eqref{w=z}}{\Rightarrow}\non\\
\pounds_{\bs{\eta}_\kappa}\bs{\sigma}^\rho=\tensor{C}{^\rho_{\kappa\tau}}\bs{\sigma}^\tau\label{lieησ}.
\end{gather}

From the above equation we can compute the exterior derivative $\bs{d\sigma}^\kappa$ of the basis $1$-forms $\{\bs{\sigma}^\kappa\}$, with the help of Cartan's identity $\pounds_{\bs{\eta}}\bs{\sigma}=\bs{d}\left(\bs{\eta}\con\bs{\sigma}\right) +\bs{\eta}\con\bs{d\sigma}$ and the fact that $\bs{\eta}\con\left(\bs{\sigma}\wedge\bs{\omega}\right)= \left(\bs{\eta}\con\bs{\sigma}\right)\wedge\bs{\omega} -\bs{\sigma}\wedge\left(\bs{\eta}\con\bs{\omega}\right)$
\bal
\pounds_{\bs{\eta}_\kappa}\bs{\sigma}^\lambda&=\tensor{C}{^\lambda_{\kappa\tau}}\bs{\sigma}^\tau\Rightarrow\non\\ \bs{d}\left(\bs{\eta}_\kappa\con\bs{\sigma}^\lambda\right)+ \bs{\eta}_\kappa\con\bs{d\sigma}^\lambda&=\tensor{C}{^\lambda_{\kappa\tau}}\bs{\sigma}^\tau\overset{\eqref{contr}}{\Rightarrow}\non\\
\bs{\eta}_\kappa\con\bs{d\sigma}^\lambda&=\tensor{C}{^\lambda_{\kappa\tau}}\bs{\sigma}^\tau\non\\
&=\tensor{C}{^\lambda_{\rho\tau}} \left(\bs{\eta}_\kappa\con\bs{\sigma}^\rho\right)\wedge\bs{\sigma}^\tau\non\\
&=\tensor{C}{^\lambda_{\rho\tau}} \left(\bs{\eta}_\kappa\con\left(\bs{\sigma}^\rho\wedge\bs{\sigma}^\tau\right) +\bs{\sigma}^\rho\wedge\left(\bs{\eta}_\kappa\con\bs{\sigma}^\tau\right)\right) \non\\
&=\tensor{C}{^\lambda_{\rho\tau}} \bs{\eta}_\kappa\con\left(\bs{\sigma}^\rho\wedge\bs{\sigma}^\tau\right) + \tensor{C}{^\lambda_{\rho\tau}}\left(\bs{\eta}_\kappa\con\bs{\sigma}^\tau\right) \wedge\bs{\sigma}^\rho\non\\
&=\tensor{C}{^\lambda_{\rho\tau}} \bs{\eta}_\kappa\con\left(\bs{\sigma}^\rho\wedge\bs{\sigma}^\tau\right)- \bs{\eta}_\kappa\con\bs{d\sigma}^\lambda\Rightarrow\non\\
2\bs{\eta}_\kappa\con\bs{d\sigma}^\lambda&=\tensor{C}{^\lambda_{\rho\tau}} \bs{\eta}_\kappa\con\left(\bs{\sigma}^\rho\wedge\bs{\sigma}^\tau\right),
\eal
and using \eqref{w=z} once more we have
\bal\label{derσ}
\bs{d\sigma}^\lambda=\frac{1}{2}\,\tensor{C}{^\lambda_{\rho\tau}} \bs{\sigma}^\rho\wedge\bs{\sigma}^\tau.
\eal

With the aid of the basis $\{\bs{\sigma}^\kappa\}$ we can define a metric tensor on $G$ by
\bal\label{metric γ}
\bs{\gamma}=\gamma_{\alpha\beta}(x^\mu)\bs{\sigma}^\alpha\otimes\bs{\sigma}^\beta,
\eal
with arbitrary functions $\gamma_{\alpha\beta}(x^\mu)$. Our concern is the \emph{homogeneous manifolds} so we demand that the manifold $G$ be \emph{invariant} under the group $G$. This means that the vector fields $\bs{\xi}_\kappa$ are also Killing fields for the metric tensor \eqref{metric γ} (of course we could choose the $\bs{\eta}$'s as Killing fields).

Applying the Killing equation to \eqref{metric γ} we have
\bga
\pounds_{\bs{\xi}_\kappa}\bs{\gamma}=0\Rightarrow \pounds_{\bs{\xi}_\kappa}\gamma_{\alpha\beta}\bs{\sigma}^\alpha\otimes \bs{\sigma}^\beta=0 \Rightarrow\non\\
\left(\pounds_{\bs{\xi}_\kappa}\gamma_{\alpha\beta}\right) \bs{\sigma}^\alpha\otimes\bs{\sigma}^\beta+ \gamma_{\alpha\beta}\left(\pounds_{\bs{\xi}_\kappa}\bs{\sigma}^\alpha\right) \otimes\bs{\sigma}^\beta+ \gamma_{\alpha\beta}\bs{\sigma}^\alpha\otimes \left(\pounds_{\bs{\xi}_\kappa}\bs{\sigma}^\beta\right)=0\overset{\eqref{lieξσ}}{\Rightarrow}\non\\
\pounds_{\bs{\xi}_\kappa}\gamma_{\alpha\beta}(x^\mu)=0\Rightarrow \tensor{\xi}{_\kappa^\lambda}\gamma_{\alpha\beta,\lambda}=0 \Rightarrow\non\\ \gamma_{\alpha\beta}(x^\mu)=const.
\ega

It is evident from the above procedure that the we use the term \emph{homogeneous manifolds} in the specific case of Lie groups with trivial isotropy subgroup.

\subsection{Constructing representations and isometries}

The local definition of an \emph{homogeneous manifold} admitting the isometries $\bs{\xi}_\kappa$ is equivalent to prescribing the metric tensor \eqref{metric γ}. The usual approach is to bring as much as possible of the $\bs{\xi}_\kappa$ into normal form (i.e. $\partial_\sigma$) and then use the first of \eqref{comξcomη} to calculate the rest of $\bs{\xi}$'s. The next step is to use \eqref{comξη} in order to find the $\bs{\eta}_\lambda$ and finally use \eqref{contr} to solve for the $\bs{\sigma}^\rho$ and thus identify the metric tensor $\bs{\gamma}$. Although the above construction is straightforward, it cannot be considered as algorithmic: it requires the solution of the system of first order partial differential equations (p.d.e.) \eqref{comξcomη} for $\bs{\xi}_\kappa$ and the subsequent action of solving the second system of p.d.e. \eqref{comξη} for $\bs{\eta}_\lambda$. It is obvious that for a large dimension $N$ of the group $G$ the above procedure can be both tedious and difficult. Below we describe a method for calculating $\bs{\xi}_\kappa$, $\bs{\eta}_\lambda$ and $\bs{\sigma}^\rho$ which does not evolve the solution of differential equations; it instead, requires only algebraic manipulations and function differentiation.

From \eqref{gen lr} it is obvious that the key element for calculating the vector fields $\bs{\xi}_\kappa,\,\bs{\eta}_\lambda$ is the knowledge of the functions $\phi^\sigma(\alpha,\beta)$. Up to now we have used the Lie group $G$ in an abstract way, thus the functions $\phi^\sigma(\alpha,\beta)$ are abstract too. In order to make the group $G$ and the functions $\phi^\sigma(\alpha,\beta)$ concrete we will use Ado's theorem \cite{Ado} which states:

\emph{Every abstract Lie algebra $\mathfrak{g}$ has a faithful representation $\rho$, on a finite-dimensional vector space $V$ of dimension $n$}.

The method we will use consists in the following steps:
\begin{enumerate}
\item With the aid of the representation $\rho$ we construct the $n\times n$ matrices $\Omega_\kappa$ corresponding to the $\bs{\xi}_\kappa$.
\item From the matrices $\Omega_\kappa$ we evaluate the finite elements of the Lie group $G$ by exponentiation
\bal\label{expO}
g_\alpha=\exp\left(\alpha^\kappa\Omega_\kappa\right)
\eal
\item The product of two elements $g_\alpha,\,g_\beta$ of $G$ returns another element $g_\phi$ of $G$, thus from the corresponding matrix multiplication
\bal\label{matrφab}
\exp\left(\phi^\gamma(\alpha,\beta)\Omega_\gamma\right)=\exp\left(\alpha^\mu\Omega_\mu\right)\times \exp\left(\beta^\nu\Omega_\nu\right)
\eal
we can read the desired functions $\phi^\gamma(\alpha,\beta)$.
\end{enumerate}

The coordinates $\alpha^\mu$ are called canonical coordinates of the first kind, or of the first genus. It is worth mentioning that we always have at our disposal the freedom to make a point transformation $\alpha^\mu=\alpha^\mu\left(\bar{\alpha}^\nu\right)$ in order to simplify the functional form of the composition function $\phi$. Later on we make repeated use of this freedom when we present the simplified form of the vector fields $\bs{\xi}_\kappa$ and $\bs{\eta}_\kappa$.

The above steps are easy to implement as soon as the representation $\rho$ is at hand. However Ado' s theorem does not give neither the representation nor the dimension $n$ of the space it acts, it only guarantees the existence of $\rho$. Thus in order to put our method at work we have to describe the way we can construct a \emph{faithful} representation $\rho$. This problem can be split in three cases depending on the properties of the Lie algebra $\mathfrak{g}$ spanned by the $\bs{\xi}_\kappa$. 

\textbf{a)} The algebra $\mathfrak{g}$ is semisimple or has a zero center; then the representation $\rho$ can be chosen as the adjoint representation, where the matrices $\Omega_\kappa$ are
\bal\label{adj}
\Omega_\kappa \in M_n(\mathbb R),\, \left(\Omega_\kappa\right)^\mu_{\ph{m}\nu}=\tensor{C}{^\mu_{\kappa\nu}}
\eal

\textbf{b)} The algebra $\mathfrak{g}$ is non-semisimple but is a direct sum of two  algebras $\mathfrak{g}=\mathfrak{h}\oplus\mathfrak{k}$; then the representation $\rho$ is the direct sum of the two representations.

{\bf c\textsubscript{1})} The algebra $\mathfrak{g}$ is non-semisimple and cannot be split into a direct sum; then by the Levi decomposition \cite{Levi} the algebra $\mathfrak{g}$ can be written as the semidirect sum $\mathfrak{g}=\mathfrak{r}\rtimes\mathfrak{h}$ where $\mathfrak{r}$ is the radical of $\mathfrak{g}$ and $\mathfrak{h}$ is a semisimple subalgebra of $\mathfrak{g}$. The desired representation $\rho$ is the extension of a faithful representation $\rho_o$ of $\mathfrak{r}$ to a representation of $\mathfrak{g}$.

{\bf c\textsubscript{2})} The algebra $\mathfrak{g}$ is a solvable and indecomposable i.e. the subalgebra $\mathfrak{h}=0$. Then we can always construct a series of subalgebras $\mathfrak{g}_1 
\subset \mathfrak{g}_2 \dots \subset \mathfrak{g}_m=\mathfrak{g}$ such that $\mathfrak{g}_{i+1}=\mathfrak{g}_i\rtimes \mathfrak{h}_i$ where $\mathfrak{h}_i$ is an one-dimensional subalgebra of $\mathfrak{g}$. The first subalgebra $\mathfrak{g}_1$ is the final algebra of the derived series of $\mathfrak{g}$ while $\mathfrak{h}_1$ is spanned by an element of the penultimate term of the derived series, and so on for the rest of $\mathfrak{h}_i$.

This method is like an upside down pyramid, since we need the knowledge of a representation of a lower dimension $\rho_o$ in order to built up $\rho$ or, to put it in other words, we need an \emph{extension} of $\rho_o$ to $\rho$. Thus we have to begin with a representation of a 2d Lie algebra to construct a representation for a 3d algebra and use it for a 4d one and so on (there are only two 2d algebras with known representations, see section \ref{BS} below). The above method serves as an alternative proof for Ado's theorem and can be found in chapter 6 of \cite{Graaf}. The basic ingredients of this method are shortly described below.

Let a Lie algebra $\mathfrak{g}=\mathfrak{s}\rtimes\mathfrak{h}$ where $\mathfrak{s}$ an ideal of $\mathfrak{g}$ and $\mathfrak{h}$ a subalgebra of $\mathfrak{g}$ and suppose we have given a faithful representation $\rho_o\colon\mathfrak{s}\to\mathfrak{gl}(V)$ of $\mathfrak{s}$ with $\text{dim}V=k$. First of all we lift up $\rho_o$ from $\mathfrak{s}$, to its universal enveloping algebra $U(\mathfrak{s})$ in the trivial way
\bal
\rho_o\left(x_1^{k_1}x_2^{k_2}\dots x_t^{k_t}\right)=\rho_o\left(x_1\right)^{k_1} \rho_o\left(x_2\right)^{k_2}\dots \rho_o\left(x_t\right)^{k_t}
\eal

Then we define an action $\Sigma\colon\mathfrak{g}\to U(\mathfrak{s})^\star$  from the algebra $\mathfrak{g}$ to the dual space of the universal enveloping algebra $U(\mathfrak{s})^\star$ of $\mathfrak{s}$ as
\bal\label{Σ action}
\Sigma\colon\mathfrak{g}\to U(\mathfrak{s})^\star,\,
\begin{cases}
\left(s\cdot f\right)(a)=f(a\,s) & s\in \mathfrak{s}\\
\left(h\cdot f\right)(a)=f(a\,h-h\,a) & h\in\mathfrak{h}
\end{cases}
\quad \forall a\in U(\mathfrak{s}),\, \forall f\in U(\mathfrak{s})^\star
\eal
and finally we consider the map
\bal
c\colon V\times V^\star\to U(\mathfrak{s})^\star,\, c\left(\upsilon,w^\star\right)(a)=w^\star\left(\rho_o(a)\upsilon\right) \quad \upsilon\in V,\,w^\star \in V^\star,\,a\in U(\mathfrak{s})
\eal

Let $\epsilon=\left\{\epsilon_1,\dots,\epsilon_k\right\}$ be a basis of $V$ and $\epsilon^\star=\left\{\epsilon_1^\star,\dots,\epsilon_k^\star\right\}$ its dual basis, e.g. $\epsilon_i^\star(\epsilon_j)=\delta_{ij}$. Then $c\left(\epsilon_i,\epsilon^\star_j\right)(a)\equiv c_{ij}(a)$ is the $(j,i)$ entry of the matrix $\rho_o(a)$, or in other words, the coefficient of $\epsilon_j$ at the expansion of the vector $\rho_o(a)\epsilon_i$ with respect to the basis $\epsilon$. Therefore each element $ c\left(\upsilon,w^\star\right)$ is called a \emph{coefficient of the representation $\rho_o$} and the image $C_{\rho_o}$ of $c$ in $U(\mathfrak{s})$ is called the \emph{coefficient space of $\rho_o$}. The representation we are seeking for, acts on a finite dimensional subspace $S_{\rho_o}\subset U(\mathfrak{s})^\star$, which is generated from a basis formed by some of the maps $c_{ij}$. This procedure is completely algorithmic and can be summarized in the following steps.
\begin{enumerate}
\item Choose a set $B=\{m_1,\dots,m_r\}$ of monomials of $U(\mathfrak{s})$ that forms a basis of the complement of $\text{ker}\rho_o$.
\item Calculate the maps $c_{ij}$, with respect to the basis $B$, which generate the coefficient space $C_{\rho_o}$.
\item Use the action \eqref{Σ action} in order to define the action of the elements of $\mathfrak{h}$ on $c_{ij}$.
\item Finally, if there are no new functions emerging from step 3 the representation space $S_{\rho_o}$ is the one spanned by $c_{ij}$, from which we can write down the desired representation. If in step 3 new coefficient functions $\xbar{c}=\{f_l,\,l\in I\}$ emerge then the representation space is $S_{\rho_o}=C_{\rho_o}\cup\{\xbar{c}\}$.
\end{enumerate}

To clarify things we give an example of the above construction for the four-dimensional Lie algebra  $A_{4,10}=\{x_1,x_2\}\rtimes\{x_3\}\rtimes\{x_4\}$ with commutation table
\bal
\left[x_2,x_3\right]=x_1,\,\left[x_3,x_4\right]=x_2,\, \left[x_2,x_4\right]=-x_3
\eal

Let $e^n_{i,j}$ denote the $n\times n$ matrices with value $1$ for the $(i,j)$ position and zero elsewhere; it is easy to see that the matrices $e^n_{1,j}$ that correspond to different values of $j$, commute with each other. Thus for an $n$-dimensional Abelian Lie algebra $\mathfrak{a}=\{x_1,\dots,x_n\}$ the matrices $\rho(x_i)=e^{n+1}_{1,i+1}$ form a representation of $\mathfrak{a}$. In our example, there is an Abelian subalgebra $\mathfrak{a}=\{x_1,x_2\}$ thus we have the representation $\rho_o(x_1)=e^3_{1,2},\,\rho_o(x_2)=e^3_{1,3}$. In order to extend $\rho_o$, firstly we choose the set $B_1=\{1,x_1,x_2\}$ of monomials of $U(\mathfrak{a})$, so the coefficients of $\rho_o$ are $c_{11},c_{12},c_{13}$ with values $c_{11}(1)=1,c_{12}(x_1)=1,c_{13}(x_2)=1$ and zero for all other values $a\in U(\mathfrak{a})$. Thus a basis for $C_{\rho_o}$ is the set $\{c_{11},c_{12},c_{13}\}$.

Let us now extend $\rho_o$ by acting with $x_3$ on $C_{\rho_o}$ with the aid of \eqref{Σ action}. Firstly, it is easy to show by induction that
\bal\label{x2lx3}
x_2^l x_3=x_3 x_2^l+l x_1 x_2^{l-1}
\eal
Thus, for an element $a=x_1^k x_2^l$ of $U(\mathfrak{a})$ and for an element $f$ of $U(\mathfrak{a})^\star$ we have
\bal\label{x3f}
x_3\cdot f(a)=f\left([a,x_3]\right)\overset{\eqref{x2lx3}}{\Rightarrow} x_3\cdot f(x_1^k x_2^l)=f\left(l x_1^{k+1} x_2^{l-1}\right)
\eal
Set $a_{kl}=l x_1^{k+1} x_2^{l-1}$.  Since for $a_{kl}$ there is not any pair $(k,l)$ resulting $a_{kl}=1$, we see from \eqref{x3f} that $x_3\cdot c_{11}=0$. As $a_{01}=x_1$ we have from \eqref{x3f} that $x_3\cdot c_{12}(x_2)=c_{12}(x_1)=1=c_{13}(x_2)$ thus $x_3\cdot c_{12}=c_{13}$. Furthermore because $a_{kl}$ never equals $x_2$ we have $x_2\cdot c_{13}=0$. Viewing the basis $\{c_{11},c_{12},c_{13}\}$ as a row vector, we have the extension of the representation $\rho_o$ to $\mathfrak{b}=\{x_1,x_2,x_3\}$
\bal
\rho_1(x_1)=e^3_{1,2},\quad \rho_1(x_2)=e^3_{1,3},\quad \rho_1(x_3)=e^3_{3,2}
\eal

Now we have to perform the extension of $\rho_1$ to our algebra $A_{4,10}$. Choosing the set $B_2=\{1,x_1,x_2,x_3\}$ of monomials of $U(\mathfrak{b})$, the coefficients of $\rho_1$ are $c_{11},c_{12},c_{13},c_{23}$ with values $c_{11}(1)=1,c_{12}(x_1)=1,c_{13}(x_2)=1,c_{23}(x_3)=1$ and zero for all other values $b\in U(\mathfrak{b})$. Thus a basis for $C_{\rho_1}$ is the set $\{c_{11},c_{12},c_{13},c_{23}\}$.

As before we now act with $x_4$ on $C_{\rho_1}$ with \eqref{Σ action} in order to extend $\rho_1$. For that we need the following relations, which again can be proven by induction
\bsub\label{rel}
\bal
x_3^m x_4&=x_4 x_3^m+mx_2 x_3^{m-1}-\frac{m(m-1)}{2}x_1 x_3^{m-2}\\
x_2^l x_4&=x_4 x_2^l-l x_2^{l-1} x_3+\frac{l(l-1)}{2}x_1 x_2^{l-2}\\
x_2^l x_3&=x_3 x_2^l+l x_1 x_2^{l-1}\\
x_3^m x_2 &=x_2 x_3^m-m x_1 x_3^{m-1}
\eal
\esub
For an element $b=x_1^k x_2^l x_3^m$ of $U(\mathfrak{b})$ and for an element $f$ of $U(\mathfrak{b})^\star$ we have
\bal\label{x4f}
x_4\cdot f(b)=f\left([b,x_4]\right)\overset{\eqref{rel}}{\Rightarrow} x_4\cdot f(x_1^k x_2^l x_3^m)=f\left(b_{klm}\right)
\eal
where
\begin{multline}\label{bklm}
b_{klm}=-l x_1^k x_2^{l-1} x_3^{m+1}+\frac{l(l-1)}{2} x_1^{k+1} x_2^{l-2} x_3^m+\\m x_1^k x_2^{l+1} x_3^{m-1}-\frac{m(m-1)}{2} x_1^{k+1} x_2^l x_3^{m-2}
\end{multline}
The element $b_{klm}$ cannot equal 1, so from \eqref{x4f} we have $x_4\cdot c_{11}=0$. As $b_{020}=x_1$  \eqref{x4f} gives $x_4\cdot c_{12}(x_2^2)=c_{12}(-2x_2 x_3+x_1)=c_{12}(-2 x_2 x_3)+c_{12}(x_1)=1$. This result forces us to define a new function $f_1\colon U(\mathfrak{b})\to\IR$ not contained in $C_{\rho_1}$ with $f_1(x_2^2)=1$; resulting $x_4\cdot c_{12}=f_1$. Likewise since $b_{002}=x_1$ we have $x_4\cdot c_{12}(x_3^2)=c_{12}(2 x_2 x_3-x_1)\Rightarrow f_1(x_3^2)=-1$. Inserting $b_{001}=x_2$ into \eqref{x4f} we have $x_4\cdot c_{13}(x_3)=c_{13}(x_2)=1$ thus $x_4\cdot c_{13}=c_{23}$, while for $b_{010}=-x_3$, $x_4\cdot c_{32}(x_2)=c_{32}(-x_3)$ resulting to  $x_4\cdot c_{32}=-c_{13}$. The element $b_{011}=-x_3^2+x_2^2$ yields to $x_4\cdot f_1(x_2 x_3)=f_1(-x_3^2+x_2^2)=-f_1(x_3^2)+f_1(x_2^2)=2$. A second function $f_2\colon U(\mathfrak{b})\to\IR$ is deduced from that relation, with $f_2(x_2 x_3)=1$ which leads to $x_4\cdot f_1=2f_2$. The action of $x_4$ on $f_2$ can be computed with the elements $b_{020}=-2 x_2 x_3+x_1$ or $b_{002}=-b_{020}$; from \eqref{x4f} we have $x_4\cdot f_2(x^2_2)=f_2(-2x_2 x_3+x_1)=-2$ yielding to  $x_4\cdot f_2=-2f_1$ (the element $b_{002}$ does not provide any further constraint). Gathering the above results we have the $\rho$-representation of $x_4$
\bal
\rho(x_4)=-e^6_{3,4}+e^6_{4,3}+e^6_{5,2}-2e^6_{5,6}+2e^6_{6,5}
\eal

To complete the representation $\rho$ we must act with each of $\{x_1,x_2,x_3\}$ on $\{f_1,f_2\}$ via \eqref{Σ action}. Firstly, $x_1\cdot f(b)=f(bx_1)$ thus $x_1\cdot f(x_1^k x_2^l x_3^m)=f\left(d_{klm}\right)$ with $d_{klm}=x_1^{k+1}x_2^lx_3^m$. Since $b_{klm}$ cannot equal either of $(x_2^2,x_3^2,x_2x_3)$ we have $x_1\cdot f_1=0$ and $x_1\cdot f_2=0$. The action of $x_2$ is $x_2\cdot f(b)=f(bx_2)$ so $x_2\cdot f(x_1^k x_2^l x_3^m)=f(e_{klm})$ where $e_{klm}$ can be computed with the aid of  \eqref{rel}; $e_{klm}=x_1^k x_2^{l+1}x_3^m-m x_1^{k+1} x_2^l x_3^{m-1}$.
For the element $b_{010}$ we have $x_2\cdot f_1(x_2)=f_1(x_2^2)=1$ which gives $x_2\cdot f_1=c_{13}$ whilst for the element $b_{001}$ we have $x_2\cdot f_2(x_3)=f_2(x_2x_3)=1$ yielding to $x_2\cdot f_2=c_{23}$. Finally, for the action of $x_3$ we have $x_3\cdot f(b)=f(bx_3)$ thus $x_3\cdot f(x_1^kx_2^l x_3^m)=f(w_{klm})$ where $w_{klm}=x_1^k x_2^l x_3^{m+1}$. Applying this action on $b_{001}$ we have $x_3\cdot f_1(x_3)=f_1(x_3^2)=-1$ which implies $x_3\cdot f_1=-c_{23}$; whilst for the element $b_{010}$ we see $x_3\cdot f_2(x_2)=f_2(x_2 x_3)=1$ thus  $x_3\cdot f_2=c_{13}$.

Finally, we end up with the set of original functions $A=\{c_{11},c_{12},c_{13},c_{23}\}$ along with the set of the new functions $\xbar{c}=\{f_1,f_2\}$. Thus, the representation space $S_{\rho_o}$ is spanned by the set $A\cup\xbar{c}$ and the desired representation is
\bal\label{A4,10rep}
\rho(x_1)&=e^6_{1,2}&\rho(x_2)&=e^6_{1,3}+e^6_{3,5}-e^6_{4,2}+e^6_{4,6}\non\\
\rho(x_3)&=e^6_{1,4} +e^6_{3,6}-e^6_{4,5}&\rho(x_4)&=-e^6_{3,4}+e^6_{4,3}+e^6_{5,2}-2e^6_{5,6}+2e^6_{6,5}
\eal
which can be used to define the desired functions $\phi^\mu(\alpha^\kappa,\beta^\kappa)$. In order to do that we have to find the finite matrix elements \eqref{expO} of the Lie group by exponentiation of the representation matrices $\Omega_\mu=\rho(x_\mu)$, i.e. $g_\alpha=\exp\left(\alpha^\mu\Omega_\mu\right)$.

This calculation can be done either by the computer algebra program "Wolfram Mathematica" or by hand, the details are given to the appendix. The result for $g_\alpha$ is 
\bal\label{matrixA410}
\begin{pmatrix}
1&\alpha^1 & \alpha^2 & \alpha^3 & \frac{1}{2}\left( -(\alpha^2)^2+(\alpha^3)^2  \right) &\alpha^2\alpha^3 \cr
0 & 1 & 0 & 0 & 0 & 0 \cr
0 & \alpha^3\cos\alpha^4 & \cos\alpha^4 & -\sin\alpha^4 &-\alpha^2\cos\alpha^4 -\alpha^3\sin\alpha^4 & -\alpha^2\sin\alpha^4+\alpha^3\cos\alpha^4\cr
0 & \alpha^3\sin\alpha^4& \sin\alpha^4 & \cos\alpha^4 & -\alpha^2\sin\alpha^4+\alpha^3\cos\alpha^4 & \alpha^2\cos\alpha^4+\alpha^3\sin\alpha^4\cr
0 & \frac{1}{2} \sin 2\alpha^4 & 0 & 0 & \cos 2\alpha^4 & \sin 2 \alpha^4\cr
0 & -\sin^2\alpha^4 & 0 & 0 & -\sin 2\alpha^4 & \cos 2\alpha^4
\end{pmatrix}
\eal

The functions \eqref{φφ σχέση} $\phi^\alpha\left(\alpha,\beta\right)$ can be found from \eqref{matrφab}
\bal
\phi^1&=\alpha^1+\beta^1+\alpha^2\beta^3\cos\beta^4+ \alpha^3\sin\beta^4\left( \beta^3-\alpha^2\sin\beta^4 \right)+ \frac{1}{4}\left(-(\alpha^2)^2+(\alpha^3)^2\right)\sin 2\beta^4\non\\
\phi^2&=\alpha^3\sin\beta^4+\alpha^2\cos\beta^4+\beta^2\non\\
\phi^3&=-\alpha^2\sin\beta^4+\alpha^3\cos\beta^4+\beta^3\non\\
\phi^4&=\alpha^4+\beta^4
\eal
The Killing fields $\bs{\xi}_\lambda$ and the invariant basis $\bs{\eta}_\mu$ can now be obtained by \eqref{gen lr}, i.e.
\bal
\left\{\bs{\xi}_\alpha\right\}&=\left\{\partial_x,\,\partial_y,\,y\partial_x+\partial_z,\,\disfrac{1}{2}(-y^2+z^2)\partial_x+z\partial_y-y\partial_z+\partial_w\right\}\non\\
\left\{\bs{\eta}_\alpha\right\}&=\left\{\partial_x,\,z\cos w \partial_x+\cos w\partial_y-\sin w\partial_z,\,z\sin w\partial_x+\sin w\partial_y+\cos w\partial_z,\,\partial_w\right\}\non
\eal
while the basis 1-forms $\bs{\sigma}^\alpha$ are calculated from \eqref{contr} 
\bal
\left\{\bs{\sigma}^\alpha\right\}=\left\{\bs{d}x-z\bs{d}y,\,\cos w\bs{d}y-\sin w\bs{d}z,\,\sin w\bs{d}y+\cos w\bs{d}z,\,\bs{d}w\right\}.
\eal

\section{Basis sets for Killing fields and 1-forms}

In this section we analyze all 2d, 3d and 4d real Lie algebras. We present the faithful representations of these algebras along with the Killing fields $\bs{\xi}_\alpha$ (left translations), the invariant fields $\bs{\eta}_\alpha$ (right translations) and the basis 1-forms $\bs{\sigma}^\alpha$. For the 3d cases we follow the Bianchi classification \cite{Bianchi} which are reproduced in \cite{Patera:invariants} and \cite{Patera:Subalg}, but we also mention other designations for these algebras, e.g. the Bianchi Type $II$ algebra is also known as Weyl algebra. For the 4d algebras we follow the conventions of \cite{Patera:invariants} and \cite{Patera:Subalg}.

The presentation is given in a form of an itemized list for each Lie algebra, with the following structure:
\begin{itemize}
\item Structure constants $\tensor{C}{^\alpha_{\beta\gamma}}$.
\item Representation $\rho(x_i)$.
\item Killing vector fields $\bs{\xi}_\alpha$.
\item Vector fields of the invariant basis $\bs{\eta}_\alpha$.
\item Basis 1-forms $\bs{\sigma}^\alpha$.
\end{itemize}

For each algebra only the nonzero structure constants $\tensor{C}{^\alpha_{\beta\gamma}}$ are given, and from them only the ones with $\beta<\gamma$.

\subsection{Two-dimensional algebras}\label{BS}

\subsubsection{$2A_1$}
\begin{itemize}\label{2A1}
\item --

\item  $\rho(x_i)=e^2_{i,i}\quad i=1,2$.

\item $\partial_x,\,\partial_y$.

\item  $\partial_x,\,\partial_y$

\item $\bs{d}x,\,\bs{d}y$.
\end{itemize}

\subsubsection{$A_2$}
\begin{itemize}\label{A2}
\item $\tensor{C}{^2_{12}}=1$.

\item $\rho(x_1)=e^2_{1,1},\,\rho(x_2)=e^2_{1,2}$.

\item $\partial_x,\,e^x\partial_y$.

\item $\partial_x+y\partial_y,\,\partial_y$.

\item $\bs{d}x,\,-y\bs{d}x+\bs{d}y$.
\end{itemize}

\subsection{Three-dimensional algebras}

\subsubsection{Bianchi Type $I$}
\begin{itemize}\label{Type I}
\item --

\item  $\rho(x_i)=e^3_{i,i}\quad i=1,2,3$.

\item $\partial_x,\,\partial_y,\,\partial_z$.

\item $\partial_x,\,\partial_y,\,\partial_z$.

\item $\bs{d}x,\,\bs{d}y,\,\bs{d}z$.
\end{itemize}

\subsubsection{Bianchi Type $II$ - $A_{3,1}$ - Weyl}
\begin{itemize}\label{Type II}
\item  $\tensor{C}{^1_{23}}=1$.

\item $\rho(x_1)=e^3_{1,2},\,\rho(x_2)=-e^3_{3,2},\,\rho(x_3)=e^3_{1,3}$.

\item $\partial_x,\,-z\partial_x+\partial_y,\,\partial_z$.

\item $\partial_x,\,\partial_y,\,-y\partial_x+\partial_z$.

\item $\bs{d}x+y\bs{d}z,\,\bs{d}y,\,\bs{d}z$.
\end{itemize}

\subsubsection{Bianchi Type $III$ - $A_1\oplus A_2$}
\begin{itemize}\label{Type III}
\item $\tensor{C}{^2_{12}}=1$.

\item $\rho(x_1)=e^3_{1,1},\,\rho(x_2)=e^3_{1,2},\,\rho(x_3)=e^3_{3,3}$.

\item $\partial_x,\,e^x\partial_y,\,\partial_z$.

\item $\partial_x+y\partial_y,\,\partial_y,\,\partial_z$.

\item $\bs{d}x,\,-y\bs{d}x+\bs{d}y,\,\bs{d}z$.
\end{itemize}

\subsubsection{Bianchi Type $IV$ - $A_{3,2}$}
\begin{itemize}\label{Type IV}
\item $\tensor{C}{^1_{13}}=1,\,\tensor{C}{^1_{23}}=1,\,\tensor{C}{^2_{23}}=1$.

\item $\rho(x_1)=e^3_{1,2},\,\rho(x_2)=e^3_{1,3},\,\rho(x_3)=e^3_{2,2}+e^3_{3,2}+e^3_{3,3}$.

\item $\partial_x,\,\partial_y,\,(x+y)\partial_x+y\partial_y+\partial_z$.

\item $e^z\partial_x,\,ze^z\partial_x+e^z\partial_y,\,\partial_z$.

\item $e^{-z}\bs{d}x-ze^{-z}\bs{d}y,\,e^{-z}\bs{d}y,\,\bs{d}z$.
\end{itemize}

\subsubsection{Bianchi Type $V$ - $A_{3,3}$ - $D\Box T_2$}
\begin{itemize}\label{Type V}
\item $\tensor{C}{^1_{13}}=1,\,\tensor{C}{^2_{23}}=1$.

\item $\rho(x_1)=e^3_{1,2},\,\rho(x_2)=e^3_{1,3},\,\rho(x_3)=e^3_{2,2}+e^3_{3,3}$.

\item $\partial_x,\,\partial_y,\,x\partial_x+y\partial_y+\partial_z$.

\item $e^z\partial_x,\,e^z\partial_y,\,\partial_z$.

\item $e^{-z}\bs{d}x,\,e^{-z}\bs{d}y,\,\bs{d}z$.
\end{itemize}

\subsubsection{Bianchi Type $VI_h$ - $A_{3,5}^h$ and $\left(h=-1\,-\,A_{3,4}
\, \text{or} \, E(1,1)\right)$}
\begin{itemize}\label{Type VI}
\item $\tensor{C}{^1_{13}}=1,\,\tensor{C}{^2_{23}}=h$.

\item $\rho(x_1)=e^3_{1,2},\,\rho(x_2)=e^3_{1,3},\,\rho(x_3)=e^3_{2,2}+h e^3_{3,3}$.

\item $\partial_x,\,\partial_y,\,x\partial_x+hy\partial_y+\partial_z$.

\item $e^z\partial_x,\,e^{h z}\partial_y,\,\partial_z$.

\item $e^{-z}\bs{d}x,\,e^{-h z}\bs{d}y,\,\bs{d}z$.
\end{itemize}

\subsubsection{Bianchi Type $VII_h$ - $A_{3,7}^h$ and $\left(h=0\,-\,A_{3,6}
\, \text{or} \, E(2)\right)$}
\begin{itemize}\label{Type VII}
\item $\tensor{C}{^1_{13}}=h,\,\tensor{C}{^2_{13}}=-1,\,\tensor{C}{^1_{23}}=1,\,\tensor{C}{^2_{23}}=h$.

\item $\rho(x_1)=e^3_{1,2},\,\rho(x_2)=e^3_{1,3},\,\rho(x_3)=h e^3_{2,2}-e^3_{2,3}+e^4_{3,2}+h e^3_{3,3}$.

\item $\partial_x,\,\partial_y,\, (h x+y)\partial_x+(-x+h y)\partial_y+\partial_z$.

\item $e^{h z}\left(\cos z\partial_x-\sin z\partial_y\right),\,e^{h z}\left(\sin z\partial_x+\cos z\partial_y\right),\,\partial_z$.

\item $e^{-h z}\left(\cos z\bs{d}x-\sin z\bs{d}y\right),\,e^{-h z}\left(\sin z\bs{d}x+\cos z\bs{d}y\right),\,\bs{d}z$.
\end{itemize}

\subsubsection{Bianchi Type $VIII$ - $A_{3,8}$}
\begin{itemize}\label{Type VIII}
\item $\tensor{C}{^1_{23}}=-1,\,\tensor{C}{^2_{13}}=-1, \,\tensor{C}{^3_{12}}=1$.

\item $\rho(x_1)=-e^3_{2,3}+e^3_{3,2},\,\rho(x_2)=-e^3_{1,3}-e^3_{3,1},\,\rho(x_3)=e^3_{1,2}+e^3_{2,1}$.

\item $\xi_{(1)},\,\xi_{(2)},\,\partial_z$.

\item $\partial_x,\,\eta_{(1)},\,\eta_{(2)}$.

\item $\bs{d}x-\sinh y\bs{d}z, \,\cos x\bs{d}y-\sin x\cosh y\bs{d}z,
\,\sin x\bs{d}y+\cos x\cosh y\bs{d}z$.
\end{itemize}

where
\bal
\bs{\xi}_{(1)}&=\sech y\cosh z\partial_x+\sinh z\partial_y-\tanh y\cosh z\partial_z \non\\
\bs{\xi}_{(2)}&=\sech y\sinh z\partial_x+\cosh z\partial_y-\tanh y\sinh z\partial_z\non\\
\bs{\eta}_{(1)}&=-\sin x\tanh y\partial_x+\cos x\partial_y-\sin x\,\sech y\partial_z\non\\
\bs{\eta}_{(2)}&=\cos x\tanh y\partial_x+\sin x\partial_y+\cos x\,\sech y\partial_z\non
\eal

\subsubsection{Bianchi Type $IX$ - $A_{3,9}$}
\begin{itemize}\label{Type IX}
\item $\tensor{C}{^1_{23}}=1,\,\tensor{C}{^2_{13}}=-1,\, \tensor{C}{^3_{12}}=1$.

\item $\rho(x_1)=-e^3_{2,3}+e^3_{3,2},\,\rho(x_2)=e^3_{1,3}-e^3_{3,1},\,\rho(x_3)=-e^3_{1,2}+e^3_{2,1}$.

\item $\xi_{(3)},\,\xi_{(4)},\,\partial_z$.

\item $\partial_x,\,\eta_{(3)},\,\eta_{(4)}$.

\item  $\bs{d}x+\sin y\bs{d}z,\, \cos x\bs{d}y-\sin x\cos y\bs{d}z,
\,\sin x\bs{d}y+\cos x\cos y\bs{d}z$.
\end{itemize}

where
\bal
\bs{\xi}_{(3)}&=\sec y\cos z\partial_x+\sin z\partial_y-\tan y\cos z\partial_z \non\\
\bs{\xi}_{(4)}&=-\sec y\sin z\partial_x+\cos z\partial_y+\tan y\sin z\partial_z\non\\
\bs{\eta}_{(3)}&=\sin x\tan y\partial_x+\cos x\partial_y-\sin x\,\sec y\partial_z\non\\
\bs{\eta}_{(4)}&=-\cos x\tan y\partial_x+\sin x\partial_y+\cos x\,\sec y\partial_z\non
\eal

\subsection{Four-dimensional algebras}
\subsubsection{$4A_1$}
\begin{itemize}\label{4A1}
\item --

\item $\rho(x_i)=e^4_{i,i}\quad i=1,\dots,4$.

\item $\partial_x,\,\partial_y,\,\partial_z,\,\partial_w$.

\item $\partial_x,\,\partial_y,\,\partial_z,\,\partial_w$.

\item $\bs{d}x,\,\bs{d}y,\,\bs{d}z,\,\bs{d}w$.
\end{itemize}

\subsubsection{$A_2\oplus 2A_1$}
\begin{itemize}\label{A22A1}
\item $\tensor{C}{^2_{12}}=1$.

\item $\rho(x_1)=e^4_{1,1},\,\rho(x_2)=e^4_{1,2},\,\rho(x_3)=e^4_{3,3},\,\rho(x_4)=e^4_{4,4}$.

\item $\partial_x,\,e^x\partial_y,\,\partial_z,\,\partial_w$.

\item $\partial_x+y\partial_y,\,\partial_y,\,\partial_z,\,\partial_w$.

\item $\bs{d}x,\,-y\bs{d}x+\bs{d}y,\,\bs{d}z,\,\bs{d}w$.
\end{itemize}

\subsubsection{$2A_2$}
\begin{itemize}\label{2A}
\item $\tensor{C}{^2_{12}}=1,\,\tensor{C}{^4_{34}}=1$.

\item $\rho(x_1)=e^4_{1,1},\,\rho(x_2)=e^4_{1,2},\,\rho(x_3)=e^4_{3,3},\,\rho(x_4)=e^4_{3,4}$.

\item $\partial_x,\,e^x\partial_y,\,\partial_z,\,e^z\partial_w$.

\item $\partial_x+y\partial_y,\,\partial_y,\,\partial_z+w\partial_w,\,\partial_w$.

\item $\bs{d}x,\,-y\bs{d}x+\bs{d}y,\,\bs{d}z,\,-w\bs{d}z+\bs{d}w$.
\end{itemize}

\subsubsection{$A_{3,1}\oplus A_1$}
\begin{itemize}\label{A31A1}
\item $\tensor{C}{^1_{23}}=1$.

\item $\rho(x_1)=e^4_{1,2},\,\rho(x_2)=-e^4_{3,2},\,\rho(x_3)=e^4_{1,3},\,\rho(x_4)=e^4_{4,4}$.

\item $\partial_x,\,-z\partial_x+\partial_y,\,\partial_z,\,\partial_w$.

\item $\partial_x,\,\partial_y,\,-y\partial_x+\partial_z,\,\partial_w$.

\item $\bs{d}x+y\bs{d}z,\,\bs{d}y,\,\bs{d}z,\,\bs{d}w$.
\end{itemize}

\subsubsection{$A_{3,2}\oplus A_1$}
\begin{itemize}\label{A32A1}
\item $\tensor{C}{^1_{13}}=1,\,\tensor{C}{^1_{23}}=1,\,\tensor{C}{^2_{23}}=1$.

\item $\rho(x_1)=e^4_{1,2},\,\rho(x_2)=e^4_{1,3},\,\rho(x_3)=e^4_{2,2}+e^4_{3,2}+e^4_{3,3},\,\rho(x_4)=e^4_{4,4}$.

\item $\partial_x,\,\partial_y,\,(x+y)\partial_x+y\partial_y+\partial_z,\,\partial_w$.

\item $e^z\partial_x,\,e^z\left(z\partial_x+\partial_y\right),\,\partial_z,\,\partial_w$.

\item $e^{-z}\left(\bs{d}x-z\bs{d}y\right),\,e^{-z}\bs{d}y,\,\bs{d}z,\,\bs{d}w$.
\end{itemize}

\subsubsection{$A_{3,3}\oplus A_1$}
\begin{itemize}\label{A33A1}
\item $\tensor{C}{^1_{13}}=1,\,\tensor{C}{^2_{23}}=1$.

\item $\rho(x_1)=e^4_{1,2},\,\rho(x_2)=e^4_{1,3},\,\rho(x_3)=e^4_{2,2}+e^4_{3,3},\,\rho(x_4)=e^4_{4,4}$.

\item $\partial_x,\,\partial_y,\,x\partial_x+y\partial_y+\partial_z,\,\partial_w$.

\item $e^z\partial_x,\,e^z\partial_y,\,\partial_z,\,\partial_w$.

\item $e^{-z}\bs{d}x,\,e^{-z}\bs{d}y,\,\bs{d}z,\,\bs{d}w$.
\end{itemize}

\subsubsection{$A_{3,4}\oplus A_1$}
\begin{itemize}\label{A34A1}
\item $\tensor{C}{^1_{13}}=1,\,\tensor{C}{^2_{23}}=-1$.

\item $\rho(x_1)=e^4_{1,2},\,\rho(x_2)=e^4_{1,3},\,\rho(x_3)=e^4_{2,2}-e^4_{3,3},\,\rho(x_4)=e^4_{4,4}$.

\item $\partial_x,\,e^z\partial_y,\,x\partial_x+\partial_z,\,\partial_w$.

\item $e^z\partial_x,\,\partial_y,\,y\partial_y+\partial_z,\,\partial_w$.

\item $e^{-z}\bs{d}x,\,\bs{d}y-y\bs{d}z,\,\bs{d}z,\,\bs{d}w$.
\end{itemize}

\subsubsection{$A^\alpha_{3,5}\oplus A_1,\, 0<|\alpha|<1$}
\begin{itemize}\label{A35A1}
\item $\tensor{C}{^1_{13}}=1,\,\tensor{C}{^2_{23}}=\alpha$.

\item $\rho(x_1)=e^4_{1,2},\,\rho(x_2)=e^4_{1,3},\,\rho(x_3)=e^4_{2,2}+\alpha e^4_{3,3},\,\rho(x_4)=e^4_{4,4}$.

\item $\partial_x,\,\partial_y,\,x\partial_x+\alpha y\partial_y+\partial_z,\,\partial_w$.

\item $e^z\partial_x,\,e^{\alpha z}\partial_y,\,\partial_z,\,\partial_w$.

\item $e^{-z}\bs{d}x,\,e^{-\alpha z}\bs{d}y,\,\bs{d}z,\,\bs{d}w$.
\end{itemize}

\subsubsection{$A_{3,6}\oplus A_1$}
\begin{itemize}\label{A36A1}
\item $\tensor{C}{^1_{23}}=1,\,\tensor{C}{^2_{13}}=-1$.

\item $\rho(x_1)=e^4_{1,2},\,\rho(x_2)=e^4_{1,3},\,\rho(x_3)=-e^4_{2,3}+e^4_{3,2},\,\rho(x_4)=e^4_{4,4}$.

\item $\partial_x,\,\partial_y,\,y\partial_x-x\partial_y+\partial_z,\,\partial_w$.

\item $\cos z\partial_x-\sin z\partial_y,\, \sin z\partial_x+\cos z\partial_y,\, \partial_z,\,\partial_w$.

\item $\cos z\bs{d}x-\sin z\bs{d}y,\,\sin z\bs{d}x+\cos z\bs{d}y,\,\bs{d}z,\,\bs{d}w$.
\end{itemize}

\subsubsection{$A^\alpha_{3,7}\oplus A_1,\, 0<\alpha$}
\begin{itemize}\label{A37A1}
\item $\tensor{C}{^1_{13}}=\alpha,\,\tensor{C}{^2_{13}}=-1,\,\tensor{C}{^1_{23}}=1,\,\tensor{C}{^2_{23}}=\alpha$.

\item $\rho(x_1)=e^4_{1,2},\,\rho(x_2)=e^4_{1,3},\,\rho(x_3)=\alpha e^4_{2,2}-e^4_{2,3}+e^4_{3,2}+\alpha e^4_{3,3},\,\rho(x_4)=e^4_{4,4}$.

\item $\partial_x,\,\partial_y,\, (\alpha x+y)\partial_x+(-x+\alpha y)\partial_y+\partial_z,\,\partial_w$.

\item $e^{\alpha z}\left(\cos z\partial_x-\sin z\partial_y\right),\,e^{\alpha z}\left(\sin z\partial_x+\cos z\partial_y\right),\,\partial_z,\,\partial_w$.

\item $e^{-\alpha z}\left(\cos z\bs{d}x-\sin z\bs{d}y\right),\,e^{-\alpha z}\left(\sin z\bs{d}x+\cos z\bs{d}y\right),\,\bs{d}z,\,\bs{d}w$.
\end{itemize}

\subsubsection{$A_{3,8}\oplus A_1$}
\begin{itemize}\label{A38A1}
\item $\tensor{C}{^1_{23}}=-1,\,\tensor{C}{^2_{13}}=-1, \,\tensor{C}{^3_{12}}=1$.

\item $ \rho(x_1)=-e^4_{2,3}+e^4_{3,2},\,\rho(x_2)=-e^4_{1,3}-e^4_{3,1},\,\rho(x_3)=e^4_{1,2}+e^4_{2,1},\,\rho(x_4)=e^4_{4,4}$.

\item $\xi_{(1)},\,\xi_{(2)},\,\partial_z,\, \partial_w$.

\item $\partial_x,\,\eta_{(1)},\,\eta_{(2)},\,\partial_w$.

\item $\bs{d}x-\sinh y\bs{d}z, \,\cos x\bs{d}y-\sin x\cosh y\bs{d}z,
\,\sin x\bs{d}y+\cos x\cosh y\bs{d}z,\,\bs{d}w$.
\end{itemize}

where
\bal
\bs{\xi}_{(1)}&=\sech y\cosh z\partial_x+\sinh z\partial_y-\tanh y\cosh z\partial_z \non\\
\bs{\xi}_{(2)}&=\sech y\sinh z\partial_x+\cosh z\partial_y-\tanh y\sinh z\partial_z\non\\
\bs{\eta}_{(1)}&=-\sin x\tanh y\partial_x+\cos x\partial_y-\sin x\,\sech y\partial_z\non\\
\bs{\eta}_{(2)}&=\cos x\tanh y\partial_x+\sin x\partial_y+\cos x\,\sech y\partial_z\non
\eal

\subsubsection{$A_{3,9}\oplus A_1$}
\begin{itemize}\label{A39A1}
\item $\tensor{C}{^1_{23}}=1,\,\tensor{C}{^2_{13}}=-1,\, \tensor{C}{^3_{12}}=1$.

\item $\rho(x_1)=-e^4_{2,3}+e^4_{3,2},\,\rho(x_2)=e^4_{1,3}-e^4_{3,1},\,\rho(x_3)=-e^4_{1,2}+e^4_{2,1},\,\rho(x_4)=e^4_{4,4}$.

\item $\xi_{(3)},\,\xi_{(4)},\,\partial_z,\, \partial_w$.

\item $\partial_x,\,\eta_{(3)},\,\eta_{(4)},\,\partial_w$.

\item $\bs{d}x+\sin y\bs{d}z,\, \cos x\bs{d}y-\sin x\cos y\bs{d}z,
\,\sin x\bs{d}y+\cos x\cos y\bs{d}z,\,\bs{d}w$.
\end{itemize}

where
\bal
\bs{\xi}_{(3)}&=\sec y\cos z\partial_x+\sin z\partial_y-\tan y\cos z\partial_z \non\\
\bs{\xi}_{(4)}&=-\sec y\sin z\partial_x+\cos z\partial_y+\tan y\sin z\partial_z\non\\
\bs{\eta}_{(3)}&=\sin x\tan y\partial_x+\cos x\partial_y-\sin x\,\sec y\partial_z\non\\
\bs{\eta}_{(4)}&=-\cos x\tan y\partial_x+\sin x\partial_y+\cos x\,\sec y\partial_z\non
\eal

\subsubsection{$A_{4,1}$}
\begin{itemize}\label{A41}
\item $\tensor{C}{^1_{24}}=1,\,\tensor{C}{^2_{34}}=1$.

\item $\rho(x_1)=e^4_{1,2},\,\rho(x_2)=e^4_{1,3},\,\rho(x_3)=e^4_{1,4},\,\rho(x_4)=e^4_{3,2}+e^4_{4,3}$.

\item $\partial_x,\,\partial_y,\,\partial_z,\,y\partial_x+z\partial_y+\partial_w$.

\item $\partial_x,\,w\partial_x+\partial_y,\,\disfrac{1}{2}w^2\partial_x+w\partial_y+\partial_z$.

\item $\bs{d}x-w\bs{d}y+\disfrac{1}{2}w^2\bs{d}z,\, \bs{d}y-w\bs{d}z,\, \bs{d}w$.
\end{itemize}

\subsubsection{$A^\alpha_{4,2},\, \alpha \notin \{0,1\}\quad \text{and}\quad A^1_{4,2}$}
\begin{itemize}\label{Aa42}
\item $\tensor{C}{^1_{14}}=\alpha,\,\tensor{C}{^2_{24}}=1,\, \tensor{C}{^2_{34}}=1,\,\tensor{C}{^3_{34}}=1$.

\item $\rho(x_1)=\alpha e^4_{1,4},\,\rho(x_2)=e^4_{2,4},\,\rho(x_3)=e^4_{2,4}+e^4_{3,4},\\\rho(x_4)=-\alpha e^4_{1,1}-e^4_{2,2}-e^4_{2,3}-e^4_{3,3}$.

\item $e^{-\alpha w}\partial_x,\,e^{-w}\partial_y,\, -e^{-w}\left(w\partial_y-\partial_z\right),\, \partial_w$.

\item $\partial_x,\,\partial_y,\,\partial_z,\, -\alpha x\partial_x-(y+z)\partial_y-z\partial_z+\partial_w$.

\item $\bs{d}x+\alpha x\bs{d}w,\,\bs{d}y+(y+z)\bs{d}w,\,\bs{d}z+z\bs{d}w,\,\bs{d}w$.
\end{itemize}

\subsubsection{$A_{4,3}$}
\begin{itemize}\label{A43}
\item $\tensor{C}{^1_{14}}=1,\,\tensor{C}{^2_{34}}=1$.

\item $\rho(x_1)=e^4_{1,2},\,\rho(x_2)=e^4_{1,3},\,\rho(x_3)=e^4_{1,4},\,\rho(x_4)=e^4_{2,2}+e^4_{4,3}$.

\item $\partial_x,\,\partial_y,\,\partial_z,\,x\partial_x+z\partial_y+\partial_w$.

\item $e^w\partial_x,\,\partial_y,\,w\partial_y+\partial_z,\,\partial_w$.

\item $e^{-w}\bs{d}x,\,\bs{d}y-w\bs{d}z,\,\bs{d}z,\,\bs{d}w$.
\end{itemize}

\subsubsection{$A_{4,4}$}
\begin{itemize}\label{A44}
\item $\tensor{C}{^1_{14}}=1,\,\tensor{C}{^1_{24}}=1,\,\tensor{C}{^2_{24}}=1,\,\tensor{C}{^2_{34}}=1,\,\tensor{C}{^3_{34}}=1$.

\item $\rho(x_1)=e^4_{1,4},\,\rho(x_2)=e^4_{1,4}+e^4_{2,4},\,\rho(x_3)=e^4_{2,4}+e^4_{3,4},\,$
\vspace{2 mm}
$\\ \rho(x_4)=-e^4_{1,1}-e^4_{1,2}-e^4_{2,2}-e^4_{2,3}-e^4_{3,3}$.

\item $e^{-w}\partial_x,\,-e^{-w}\left(w\partial_x-\partial_y\right),\, e^{-w}\left(\disfrac{1}{2}w^2\partial_x-w\partial_y+\partial_z\right),\,\partial_w$.

\item $\partial_x,\,\partial_y,\,\partial_z,\,-(x+y)\partial_x-(y+z)\partial_y-z\partial_z+\partial_w$.

\item  $\bs{d}x+(x+y)\bs{d}w,\,\bs{d}y+(y+z)\bs{d}w,\,\bs{d}z+z\bs{d}w,\,\bs{d}w$.
\end{itemize}

\subsubsection{$\left(A^{\alpha,\beta}_{4,5},\,\alpha\beta\neq0,\, -1\leq\alpha<\beta<1\right) \, \text{and} \, \left(A^{\alpha,\alpha}_{4,5},\,\alpha\neq0,\, -1\leq\alpha<1\right) \, \text{and} \\ \left(A^{\alpha,1}_{4,5},\,\alpha\neq0,\, -1\leq\alpha<1\right)  \, \text{and} \, \left(A^{1,1}_{4,5},\,\alpha\neq0\right)$}

\begin{itemize}\label{Aab45}
\item  $\tensor{C}{^1_{14}}=1,\,\tensor{C}{^2_{24}}=\alpha,\,\tensor{C}{^3_{34}}=\beta$.

\item $ \rho(x_1)=e^4_{1,4},\,\rho(x_2)=\alpha e^4_{2,4},\,\rho(x_3)=\beta e^4_{3,4},\,\rho(x_4)=-e^4_{1,1}-\alpha e^4_{2,2}-\beta e^4_{3,3}$.

\item $e^{-w}\partial_x,\,e^{-a w}\partial_y,\,e^{-\beta w}\partial_z,\,\partial_w$.

\item $\partial_x,\,\partial_y,\,\partial_z,\,-x\partial_x-\alpha y\partial_y-\beta z\partial_z+\partial_w$.

\item $\bs{d}x+x\bs{d}w,\,\bs{d}y+\alpha y\bs{d}w,\,\bs{d}z+\beta z\bs{d}w,\,\bs{d}w$.
\end{itemize}

\subsubsection{$A^{\alpha,\beta}_{4,6},\,\alpha\neq 0,\,\beta\geq0$}
\begin{itemize}\label{Aab46}
\item $\tensor{C}{^1_{14}}=\alpha,\,\tensor{C}{^2_{24}}=\beta,\,\tensor{C}{^2_{34}}=1,\,\tensor{C}{^3_{24}}=-1,\,\tensor{C}{^3_{34}}=\beta$.

\item $\rho(x_1)=\alpha e^4_{1,4},\,\rho(x_2)=\beta e^4_{2,4}-e^4_{3,4},\,\rho(x_3)=e^4_{2,4}+\beta e^4_{3,4},\,$
\vspace{2 mm}
$\\ \rho(x_4)=-\alpha e^4_{1,1}-\beta e^4_{2,2}+e^4_{3,2}-e^4_{2,3}-\beta e^4_{3,3}$.

\item  $e^{-\alpha w}\partial_x,\,e^{-\beta w}\left(\cos w\partial_y+\sin w\partial_z\right),\,-e^{-\beta w}\left(\sin w\partial_y-\cos w\partial_z\right),\,\partial_w$.

\item $\partial_x,\,\partial_y,\,\partial_z,\,-\alpha x\partial_x-(\beta y+z)\partial_y+(y-\beta z)\partial_z+\partial_w$.

\item $\bs{d}x+\alpha x\bs{d}w,\,\bs{d}y+(\beta y+z)\bs{d}w,\,\bs{d}z-(y-\beta z)\bs{d}w,\, \bs{d}w$.
\end{itemize}

\subsubsection{$A_{4,7}$}
\begin{itemize}\label{A47}
\item $\tensor{C}{^1_{14}}=2,\,\tensor{C}{^2_{24}}=1,\,\tensor{C}{^1_{23}}=1,\,\tensor{C}{^2_{34}}=1,\,\tensor{C}{^3_{34}}=1$.

\item $\rho(x_1)=2e^4_{1,4},\,\rho(x_2)=e^4_{1,3}+e^4_{2,4},\,\rho(x_3)=-e^4_{1,2}+e^4_{2,4}+e^4_{3,4},\,$
\vspace{2 mm}
$\\  \rho(x_4)=-2e^4_{1,1}-e^4_{2,2}-e^4_{2,3}-e^4_{3,3}$.

\item $e^{-2 w}\partial_x,\,-e^{-w}\left(\disfrac{1}{2}z\partial_x-\partial_y\right),
\, e^{-w}\left(\disfrac{1}{2}(y-z-zw)\partial_x-w\partial_y+\partial_z\right),\,\partial_w$.

\item $\partial_x,\,\disfrac{1}{2}z\partial_x+\partial_y,\,-\frac{1}{2} (y+z)\partial_x+\partial_z,\,-2x\partial_x-(y+z)\partial_y-z\partial_z+\partial_w$.

\item $\bs{d}x-\disfrac{z}{2} \bs{d}y+\disfrac{y+z}{2}\bs{d}z+2x\bs{d}w,\,\bs{d}y+(y+z)\bs{d}w,\,\bs{d}z+z\bs{d}w,\,\bs{d}w$.
\end{itemize}

\subsubsection{$A_{4,8}$}
\begin{itemize}\label{A48}
\item $\tensor{C}{^1_{23}}=1,\,\tensor{C}{^2_{24}}=1,\,\tensor{C}{^3_{34}}=-1$.

\item $\rho(x_1)=e^4_{2,3},\,\rho(x_2)=e^4_{2,4},\,\rho(x_3)=e^4_{4,3},\,\rho(x_4)=e^4_{4,4}$.

\item $\partial_x,\,\partial_y,\,y\partial_x+e^w\partial_z,\,y\partial_y+\partial_w$.

\item $\partial_x,\,z\partial_x+e^w\partial_y,\,\partial_z,\,z\partial_z+\partial_w$.

\item $\bs{d}x-e^{-w}z\bs{d}y,\,e^{-w}\bs{d}y,\,\bs{d}z-z\bs{d}w,\,\bs{d}w$.
\end{itemize}

\subsubsection{$\left(A^\beta_{4,9},\,0<|\beta|<1\right)\,\text{and}\, A^0_{4,9}\,\text{and}\, A^1_{4,9}$}
\begin{itemize}\label{Ab49}
\item $\tensor{C}{^1_{23}}=1,\,\tensor{C}{^1_{14}}=1+\beta,\,\tensor{C}{^2_{24}}=1,\,\tensor{C}{^3_{34}}=\beta$.

\item $\rho(x_1)=(1+\beta)e^4_{1,4},\,\rho(x_2)=e^4_{1,3}+e^4_{2,4},\,\rho(x_3)=-e^4_{1,2}+\beta e^4_{3,4},$
\vspace{2 mm}
$\\  \rho(x_4)=-(1+\beta)e^4_{1,1}-e^4_{2,2}-\beta e^4_{3,3}$.

\item  $e^{-(1+\beta) w}\partial_x,-e^{-w}\left(\disfrac{z}{1+\beta}\partial_x-\partial_y\right),\, e^{-\beta w}\left(\disfrac{\beta y}{1+\beta}\partial_x+\partial_z\right),\,\partial_w$.

\item $\partial_x,\,\disfrac{\beta z}{1+\beta}\partial_x+\partial_y,-\disfrac{y}{1+\beta} \partial_x+\partial_z,\, -(1+\beta)x\partial_x-y\partial_y-\beta z\partial_z+\partial_w$.

\item  $\bs{d}x-\disfrac{\beta z}{1+\beta}\bs{d}y+\disfrac{y}{1+\beta}\bs{d}z+(1+\beta)x\bs{d}w,\,\bs{d}y+y\bs{d}w,\,\bs{d}z+\beta z\bs{d}w,\,\bs{d}w$.
\end{itemize}

\subsubsection{$A_{4,10}$}
\begin{itemize}\label{A410}
\item $\tensor{C}{^1_{23}}=1,\,\tensor{C}{^2_{34}}=1,\,\tensor{C}{^3_{24}}=-1$.

\item see \eqref{A4,10rep}.

\item  $\partial_x,\,\partial_y,\,y\partial_x+\partial_z,\,\disfrac{1}{2}(-y^2+z^2)\partial_x+z\partial_y-y\partial_z+\partial_w$.

\item  $\partial_x,\,z\cos w \partial_x+\cos w\partial_y-\sin w\partial_z,\,z\sin w\partial_x+\sin w\partial_y+\cos w\partial_z,\,\partial_w$.

\item $\bs{d}x-z\bs{d}y,\,\cos w\bs{d}y-\sin w\bs{d}z,\,\sin w\bs{d}y+\cos w\bs{d}z,\,\bs{d}w$.
\end{itemize}

\subsubsection{$A_{4,11},\, \alpha>0$}
\begin{itemize}\label{A411}
\item  $\tensor{C}{^1_{23}}=1,\,\tensor{C}{^1_{14}}=2\alpha,\,\tensor{C}{^2_{24}}=\alpha,\,\tensor{C}{^3_{24}}=-1,\,\tensor{C}{^2_{34}}=1,\,\tensor{C}{^3_{34}}=\alpha$.

\item$\rho(x_1)=2\alpha e^4_{1,4},\,\rho(x_2)=e^4_{1,3}+\alpha e^4_{2,4}-e^4_{3,4},$
\vspace{2 mm}
$\\ \rho(x_3)=-e^4_{1,2}+e^4_{2,4}+\alpha e^4_{3,4},$
\vspace{2 mm}
$\\  \rho(x_4)=-2\alpha e^4_{1,1}-\alpha e^4_{2,2}-e^4_{2,3}+e^4_{3,2}-\alpha e^4_{3,3}$.

\item $\partial_x,\,\partial_y,\,(y-\disfrac{1}{2\alpha}z)\partial_x+\partial_z,\,\xi_{(5)}$.

\item $e^{2\alpha w}\partial_x,\, \eta_{(5)},\,\eta_{(6)},\partial_w$.

\item $\bs{\sigma}^1,\,\bs{\sigma}^2,\,\bs{\sigma}^3,\,\bs{d}w$.
\end{itemize}
where
\bal
\bs{\xi}_{(5)}&=\disfrac{1}{2\alpha}(-\alpha y^2+\alpha z^2+y z+4\alpha^2 x)\partial_x+(\alpha y+z)\partial_y-(y-\alpha z)\partial_z+\partial_w\non\\
\bs{\eta}_{(5)}&=e^{\alpha w}\left(\disfrac{1}{2\alpha}z\left(2\alpha\cos w+\sin w\right) \partial_x+\cos w\partial_y-\sin w\partial_z\right)\non\\
\bs{\eta}_{(6)}&=e^{\alpha w}\left(\disfrac{1}{2\alpha}z\left(-\cos w+2\alpha\sin w\right) \partial_x+\sin w\partial_y+\cos w\partial_z\right)\non\\
\bs{\sigma}^1&=e^{-2\alpha w}\left(\bs{d}x-z\bs{d}y+\frac{1}{2\alpha}z\bs{d}z\right)\non\\
\bs{\sigma}^2&=e^{-\alpha w}\left(\cos w\bs{d}y-\sin w\bs{d}z\right)\non\\
\bs{\sigma}^3&=e^{-\alpha w}\left(\sin w\bs{d}y+\cos w\bs{d}z\right)\non
\eal

\subsubsection{$A_{4,12}$}
\begin{itemize}\label{A412}
\item  $\tensor{C}{^1_{13}}=1,\,\tensor{C}{^1_{24}}=1,\,\tensor{C}{^2_{14}}=-1,\,\tensor{C}{^2_{23}}=1$.

\item $\rho(x_1)=e^4_{1,3}-e^4_{2,4},\,\rho(x_2)=e^4_{1,4}+e^4_{2,3},\,\rho(x_3)=-e^4_{1,1}-e^4_{1,2},\rho(x_4)=-e^4_{1,2}+e^4_{2,1}$.

\item  $e^{-z}\left(\cos w\partial_x+\sin w\partial_y\right),\,e^{-z}\left(-\sin w\partial_x+\cos w\partial_y\right),\,\partial_y,\,\partial_w$.

\item $\partial_x,\,\partial_y,\,-x \partial_x-y\partial_y+\partial_z,\,-y\partial_x+x\partial_y+\partial_w$.

\item  $\bs{d}x+x\bs{d}z+y\bs{d}w,\, \bs{d}y+y\bs{d}z-x\bs{d}w,\,\bs{d}z,\,\bs{d}w$.
\end{itemize}

\section{Discussion}

In this work we presented a method for constructing faithful representations of Lie algebras and showed how to use them in order to define the homogeneous manifolds via the basis 1-forms. Furthermore we used these representations for calculating the Killing and the invariant fields of the above manifolds, by identification with left and right translations respectively and the use of the composition function of the Lie group. The key difference from other works in the subject is that we \emph{do not} solve any differential equations in order to determinate the Killing/invariant fields or the 1-forms for each and every algebra; we only apply algebraic manipulations and function differentiation.

The Killing fields, $\bs{\xi}_\alpha$ can be found in the exhaustive and elegant work of Popovych et al in \cite{Popovych:RealLieAlg}; where they use the notion of megaideals to simplify the calculations, but at the end of the day they do need to solve partial differential equations. Furthermore we must mention the work of Shirokov \cite{ShirI}, \cite{ShirII} where is presented a method to calculate the Killing fields $\bs{\xi}_\alpha$. A way to evaluate the composition function $\phi$ along with the invariant vector fields $\bs{\eta}_\kappa$ was presented in \cite{Magazev} using quadratures.

The basis 1-forms $\bs{\sigma}^\alpha$ for the 4d case, were presented in \cite{Mojaveri:4+1hsc} and \cite{Jafarizadeh:PoisonLie} with the help of  Maurer--Cartan form of the corresponding Lie group, (see \cite{Olver:FramesI} for details). In these works the authors used the adjoint representation in order to calculate the 1-forms $\bs{\sigma}^\alpha$ whereas we base our presentation on a faithful representation, which of course coincides with the adjoint whenever the algebra is semisimple.

The calculations are made with the aid of computer algebra system "Wolfram Mathematica", but as is pointed out in the main text and in the appendix, these calculations can also be done by hand.

\section*{Acknowledgments}
I would to thank Prof. T. Christodoulakis for reading the manuscript and providing useful comments.

\appendix
\section{Calculation of Matrix Exponent}

The classical method for calculating functions of matrices can be found in any standard textbook on matrices, see for example \cite{Gantmacher} page 97,

Let $f(A)$ be a function defined on the spectrum of a matrix $A$ and $r(A)$ the corresponding Lagrange--Sylvester interpolation polynomial. Then $f(A) = r(A)$.

Let $r(\lambda)=\sum_{k=0}^{n-1} c_k \lambda^k$ and $\lambda_i$ be an eigenvalue of the matrix $A$; if $\lambda_i$ is a simple root of $r(\lambda)$ then $\exp\lambda_i=r(\lambda_i)$; if the multiplicity of $\lambda_i$ is $m_i$, then
\bal
e^{\lambda_i}=\frac{d^s r(\lambda)}{d\lambda^s}\Big|_{\lambda=\lambda_i} \quad s=1,2,\dots,m_i-1.
\eal

Applying the above relation to the eigenvalues of the matrix $\exp A$ we add up with a system of $n$ equations for the $n$ unknowns $c_k$.

As an example we calculate the matrix $S(\theta^\gamma)=\exp\left(\theta^\mu\Omega_\mu \right)$ where $\Omega_\mu$ are the matrices \eqref{A4,10rep} from the representation of the algebra $A_{4,10}$.

The eigenvalues $\lambda_i$, along with their multiplicities $m_i$, of the matrix $S(\theta^\gamma)$ are
\bal
\lambda_1=0,\, m_1=2 \quad \lambda_{2,3}=\pm\ima\, \theta^4,\, m_{2,3}=1 \quad \lambda_{4,5}=\pm 2\ima\, \theta^4,\, m_{3,4}=1
\eal
and the linear system that the constants obey is
\bal
1&=r(0)\Rightarrow 1=c_0\non\\
1&=r'(0)\Rightarrow 1=c_1\non\\
\exp\left(\pm\ima\,\theta^4\right)&=r(\pm\ima\,\theta^4)\Rightarrow \exp\left(\pm\ima\,\theta^4\right)= \sum_{k=0}^{5}c_k\left(\ima\,\theta^4\right)^k\non\\
\exp\left(\pm2\ima\,\theta^4\right)&=r(\pm2\ima\,\theta^4)\Rightarrow \exp\left(\pm2\ima\,\theta^4\right)= \sum_{k=0}^{5}c_k\left(\pm2\ima\,\theta^4\right)^k
\eal
The solution of the above system is
\bal
\begin{split}
c_0=1,\,c_1=1,\,c_2=\frac{1}{3(\theta^4)^2}\left(7-\cos\theta^4\right)\sin^2\frac{\theta^4}{2},\\
c_3=\frac{1}{24(\theta^4)^3}\left( 30\theta^4-32\sin\theta^4+\sin2\theta^4 \right),\, c_4=\frac{2}{3(\theta^4)^4}\sin^4\frac{\theta^4}{2}\\
c_5=\frac{1}{24(\theta^4)^5}\left( 6\theta^4-8\sin\theta^4+\sin 2\theta^4 \right)
\end{split}
\eal
With the aid of the above values we can calculate $\exp S$ from the interpolation polynomial $r(A)$. In order to simplify the resulting expression (see the paragraph after \eqref{matrφab}) we make the redefinitions
\bal
\theta^1&\mapsto \frac{1}{8}\csc^2\frac{\alpha^4}{2}\Big(  4\alpha^1-2\alpha^2\alpha^3-\left((\alpha^2)^2+(\alpha^3)^2\right)\alpha^4 \non\\
& \ph{\mapsto}- \left(4\alpha^1-2\alpha^2\alpha^3\right)\cos\alpha^4   +\left((\alpha^2)^2+(\alpha^3)^2\right)\sin\alpha^4 \Big)\non\\
\theta^2&\mapsto \frac{1}{2}\left( \alpha^2\alpha^4\cot\frac{\alpha^4}{2}-\alpha^3\alpha^4 \right)\non\\
\theta^3&\mapsto \frac{1}{2}\left( \alpha^2+\alpha^3\cot\frac{\alpha^4}{2} \right)\alpha^4
\eal
thus we arrive at the matrix \eqref{matrixA410}.

\phantomsection
\addcontentsline{toc}{section}{References}

\bibliographystyle{utphys}
\bibliography{Main}

\end{document}